\documentclass[12pt]{article}
\usepackage{amsmath}
\usepackage{latexsym}
\usepackage{amssymb}
\usepackage{a4wide}
\newtheorem{thm}{Theorem}

\newtheorem{Defn}{Definition}
\newtheorem{Remark}{Remark}
\newtheorem{prop}{Proposition}

\newenvironment{rem}{\begin{Remark}\rm}{\end{Remark}}
\newenvironment{proof}{{\noindent\bf Proof.}}%
                  {\nopagebreak\hspace*{\fill}$\Box$\medskip\medskip\par}

\newcommand{\wt}{\widetilde}

\newcommand{\tensor}{\otimes}

\newcommand{\n}{\rm}

\newcommand{\mto}{\mapsto}

\newcommand{\ve}{\varepsilon}

\newcommand{\isom}{\cong}

\newcommand{\N}{{\mathbb N}}

\newcommand{\R}{{\mathbb R}}

\newcommand{\Q}{{\mathbb Q}}

\newcommand{\Z}{{\mathbb Z}}
\newcommand{\C}{{\mathbb C}}

\newcommand{\K}{{\mathbb K}}

\newcommand{\cO}{{\cal O}}

\newcommand{\sub}{\subseteq}

\newcommand{\GL}{\mbox{\rm GL}}

\newcommand{\cB}{{\cal B}}
\newcommand{\cA}{{\mathcal A}}

\newcommand{\sbull}{{\scriptscriptstyle \bullet}}

\begin{document}
\begin{center}
{\Large \bf Tensor products in the category of topological\vspace{1mm}
vector spaces
are not associative}\vspace{3 mm}\\
{\bf Helge Gl\"{o}ckner}
\end{center}
{\bf Abstract.}
We show by example that the associative law does not hold
for tensor products in the category
of general (not necessarily locally convex)
topological vector spaces.
The same pathology occurs for
tensor products of Hausdorff abelian topological
groups.\renewcommand{\thefootnote}{\fnsymbol{footnote}}
\footnotetext{\hspace{-2mm}{\em Classification\/}: 46A16, 46A32; 22A05}
\renewcommand{\thefootnote}{\arabic{footnote}}
%
%CLASSIFIKATIONEN BEDEUTEN:
%not lcx spaces;
%top tens prod;
%structure of gen top gps
%
%
\section*{Introduction}
Let $\cA$ be a class of (not necessarily Hausdorff)
real topological vector spaces
(resp., of abelian topological groups),
such that $\cA$ is closed under the formation of
finite cartesian products.
Given $E_1,\ldots, E_n\in \cA$ (where $2\leq n\in \N$),
we call an element $T\in \cA$,
together with an $n$-linear (resp., $n$-additive)
continuous map
\[
\tau\!: E_1\times \cdots\times E_n\to T\,,
\]
a {\em tensor product
of $E_1,\ldots, E_n$ in the class~$\cA$\/} if
for every $E\in \cA$
and $n$-linear (resp., $n$-additive)
continuous map $f\!: E_1\times\cdots \times E_n\to E$,
there exists a unique continuous linear map
(resp., continuous homomorphism)
$\wt{f}\!: T\to E$ such that $\wt{f}\circ \tau=f$.\\[3mm]
For example, the tensor products in the
class of Hausdorff locally convex spaces
are the
{\em projective tensor products},
going back to
Grothendieck's memoir~\cite{Gro}.
In this case, an explicit description
of the locally convex topology
(by means of suitable cross-seminorms)
is available, and it is well-known that
an associative law holds for iterated
projective tensor products;
this is important for applications
in topological algebra.\\[3mm]
Two-fold tensor products $E\tensor F$
in the class of real topological vector spaces
have been studied in \cite{Tom}, \cite{Wae},
\cite{Iya}, \cite{Hol}, \cite{Des}
and the breakthrough
papers \cite{Tu1}, \cite{Tu2}.
For two-fold tensor products
in various classes
of abelian topological groups, see
\cite{Hof}, \cite{Gar}, \cite{Tom},
\cite{Akm}, \cite{Lot}, and \cite{Sah}.
In none of these works, higher tensor products
or iterated tensor products are discussed,
and accordingly the question of associativity of
tensor products has not been raised there.\\[3mm]
The present paper intends to close this gap.
Based on an explicit
description of the topology on tensor products
in the category of real topological vector spaces
provided in Section~1
(Proposition~\ref{propo}),
we first describe a sufficient condition
ensuring that
$(E_1\tensor E_2)\tensor E_3$
be canonically isomorphic to
$E_1\tensor (E_2\tensor E_3)$:
it suffices that the outer factors $E_1$ and $E_3$
be locally bounded (Proposition~\ref{criter}).
We then establish the main result:
For $E:=\R^\N$, none of the tensor products
$(E\tensor E)\tensor E$,
$E\tensor (E\tensor E)$
and $E\tensor E\tensor E$
(in the class of real topological vector spaces)
are naturally isomorphic (Theorem~\ref{main}).
Likewise, the associative law fails for
tensor products in the categories of
Hausdorff real topological vector spaces
and Hausdorff
abelian topological groups (Remark~\ref{conseq}).
\section{Description of the topology on tensor products}
In the following,
topological vector spaces are {\em not\/}
presumed Hausdorff, nor are topological groups,
unless we explicitly say the contrary.
\begin{prop}\label{propo}
Given real topological vector spaces
$E_1,\ldots, E_n$,
let $T:=E_1\tensor_\R \cdots \tensor_\R E_n$
be the $($algebraic$)$
tensor product of the vector spaces
$E_1,\ldots, E_n$, and
$\tau\!: E_1\times \cdots\times E_n\to T$,
$\tau(x_1,\ldots,x_n):=x_1\tensor\cdots\tensor x_n$.
Let $\cB$ be the set of all subsets~$U$
of $\,T$ of the form
\begin{equation}\label{specialneigh}
U\;=\;\sum_{k\in \N} \tau(U_{k,1}\times\cdots\times U_{k,n})
\;:=\;\bigcup_{k\in\N}\,
\sum_{\ell=1}^k \tau(U_{\ell,1}\times\cdots\times U_{\ell,n})\,,
\end{equation}
where $(U_{k,j})_{k\in \N}$ is a sequence
of balanced zero-neighbourhoods
in $E_j$, for $j=1,\ldots, n$.
Then the following holds:
\begin{itemize}
\item[\rm (a)]
$\cB$ is a basis for the filter of zero-neighbourhoods
of some topology $\cO$ on~$T$ making $T$ a real topological
vector space.
\item[\rm (b)]
$\tau\!: E_1\times\cdots\times E_n\to (T,\cO)$ is continuous.
\item[\rm (c)]
$(T,\cO)$, together with the continuous
$n$-linear map~$\tau$, is a tensor product
of $E_1,\ldots, E_n$ in the category of
real topological vector spaces.
\item[\rm (d)]
If $E_1,\ldots, E_n$ are Hausdorff,
then $T$ is Hausdorff and hence is the tensor
product of $E_1,\ldots, E_n$
in the category of Hausdorff real topological vector spaces.
\item[\rm (e)]
If $E_1,\ldots, E_n$
are Hausdorff, then $(T,\cO)$, considered as an abelian
topological group,
together with the continuous $n$-additive map $\tau$,
is a tensor product
of the Hausdorff abelian topological groups
$(E_1,+),\ldots, (E_n,+)$
in the category of Hausdorff abelian topological groups.
\end{itemize}
\end{prop}
\begin{proof}
(a) It is obvious
that every $U\in \cB$
is balanced, absorbing, and that $tU\in \cB$
for any $t\in \R^\times$;
hence
conditions
(EV${}_I$)
and (EV${}_{II}$) of
\cite[I, \S1, No.\,5, Prop.\,4]{BTV}
are satisfied.
In order that $\cB$ be a basis
of a vector topology, it remains to verify
condition
(EV${}_{III}$). To this end, let $U$ be as in (\ref{specialneigh}).
We find
balanced zero-neighbourhoods $V_{k,j}$ in $E_j$
for each $k\in \N$ and $j\in \{1,\ldots n\}$
such that
$V_{k,j}\sub U_{2k-1,j}\cap U_{2k,j}$.
Then, re-ordering terms
and abbreviating $V_{k,1}\tensor\cdots\tensor V_{k,n}:=
\tau(V_{k,1}\times\cdots\times V_{k,n})$,
we find that
\begin{eqnarray*}
\!\!\!\!\lefteqn{\left(
\sum_{k\in \N} V_{k,1}\tensor \cdots\tensor V_{k,n}\right)
\;+\;\left(\sum_{k\in \N} V_{k,1}\tensor \cdots\tensor  V_{k,n}\right)}\qquad\\
\!\!\!\!&=&
V_{1,1}\!\tensor\! \cdots\tensor V_{1,n}\,+\,
V_{1,1}\!\tensor\cdots\tensor\! V_{1,n}
\,+\,
V_{2,1}\!\tensor \cdots\tensor\! V_{2,n}\,+\,
V_{2,1}\!\tensor \cdots\tensor \!V_{2,n}+\cdots\\
&\sub & \sum_{k\in \N} U_{k,1}\tensor \cdots\tensor U_{k,n}\,.
\end{eqnarray*}
Thus $\cB$ is a basis for the filter of zero-neighbourhoods
of a vector topology~$\cO$ on~$T$.\vspace{3mm}

(b) In order that the $n$-linear map $\tau$ be
continuous, we only need to show it is continuous
at zero (see \cite{BTV}, Chap.\,I, \S1, No.\,6, Prop.\,5).
Now, given any zero neighbourhood $U\in \cB$ as in (\ref{specialneigh}),
we have
$\tau^{-1}(U)\supseteq U_{1,1}\times\cdots\times U_{1,n}$,
which is a neighbourhood of $(0,\ldots, 0)$.\vspace{3mm}

(c) Suppose that $f\!: E_1\times\cdots\times E_n\to E$ is a continuous
$n$-linear map to a real topological vector space~$E$.
Since $(T,\tau)$ is the (algebraic) tensor product
of the vector spaces $E_1,\ldots, E_n$, there is a unique
linear map $\wt{f}\!: T\to E$ such that
$\wt{f}\circ\tau=f$.
It only remains to show that $\wt{f}$ is continuous.
To this end, let $W\sub E$ be a zero-neighbourhood.
Then standard arguments provide
a sequence of zero-neighbourhoods $W_k\sub E$
such that
$\sum_{k\in \N} W_k\sub W$.
The $n$-linear map~$f$ being continuous,
for each $k\in \N$ we find balanced
zero-neighbourhoods
$U_{k,j}\sub E_j$ for $j=1,\ldots, n$
such that
$f(U_{k,1}\times\cdots\times U_{k,n})\sub W_k$.
Then $U:=\sum_{k\in \N}\tau(U_{k,1}\times\cdots\times
U_{k,n})\in \cB$ is a zero-neighbourhood in~$T$, and
\[
\wt{f}(U)=\sum_{k\in \N} \wt{f}\bigl(
\tau
(U_{k,1}\times\cdots\times U_{k,n})\bigr)
=\sum_{k\in \N} f(U_{k,1}\times\cdots\times U_{k,n})
\sub \sum_{k\in \N} W_k\sub W\,.
\]
Thus the linear map
$\wt{f}$ is continuous at zero and thus continuous.\vspace{2mm}

(d) The proof is by induction. For $n=2$,
the assertion is Turpin's celebrated result (see \cite{Tu1},
\cite{Tu2}). Now, by induction $E_1\tensor\cdots \tensor E_{n-1}$
is Hausdorff, and hence so is $F:=(E_1\tensor \cdots \tensor E_{n-1})\tensor
E_n$. The continuous $n$-linear map $f\!: E_1\times\cdots\times E_n\to F$,
$(x_1,\ldots, x_n)\mto (x_1\tensor\cdots\tensor x_{n-1})\tensor x_n$
induces a continuous linear map
$\tilde{f}\!: T=E_1\tensor\cdots \tensor E_n\to F$,
determined by $\tilde{f}\circ \tau=f$.
It is known from abstract algebra
that $\tilde{f}$ is an isomorphism of vector spaces.
Since $F$ is Hausdorff and $\tilde{f}$ is a continuous injection,
$T$ is Hausdorff.\vspace{2mm}

(e) To outline the idea,
let us assume that $n=2$ (the general case being
analogous).
Suppose that $f\!: E_1\times E_2\to A$ is
a continuous
bi-additive map to an abelian topological group~$A$.
Using the given scalar multiplication
also on the right, we consider $E_1$ as a $(\Z,\R)$-bimodule.
Similarly, $E_2$ is considered as a $(\R,\Z)$-bimodule.
Then $f$ is $\R$-balanced in the sense of \cite[p.\,161]{Coh}.
In fact, for any $x\in E_1$, $y\in E_2$
and $q\in \Q$, say $q=\frac{m}{n}$
with $m\in \Z$ and $n\in \N\setminus\{0\}$,
we have
$f(qx,y)=f\bigl(m{\textstyle \frac{1}{n}}x,n{\textstyle \frac{1}{n}}y\bigr)
=f\bigl(n{\textstyle \frac{1}{n}}x,m{\textstyle \frac{1}{n}}y\bigr)
=f(x,qy)$,
whence $f(rx,y)=f(x,ry)$ for all
$r\in \R$, by continuity.
Now, $f$ being $\R$-balanced,
there is a uniquely
determined homomorphism of groups
$\wt{f}\!: T\to A$ such that
$\wt{f}\circ\tau=f$ (cf.\ \cite[pp.\,161--162]{Coh}).
As in (c), we see that $\wt{f}$
is continuous.
\end{proof}
\section{A criterion for associativity
of tensor products}\label{seccrit}
Given
real topological vector spaces
$E_1$, $E_2$, $E_3$,
there is a continuous linear map
$\phi\!: 
E_1\tensor E_2\tensor E_3\to (E_1\tensor E_2)\tensor E_3$,
uniquely determined by $\phi(x\tensor y\tensor z)=(x\tensor y)\tensor
z$, and $\phi$ is an isomorphism of vector spaces
(cf.\ proof of Proposition~\ref{propo}\,(d)).
Likewise, there is a unique continuous
linear map (and isomorphism of vector spaces)
$\psi\!:
E_1\tensor E_2\tensor E_3\to E_1\tensor (E_2\tensor E_3)$,
determined by $\psi(x\tensor y\tensor z)=x\tensor (y\tensor z)$.\\[3mm]
The following proposition describes criteria
ensuring that $\phi$, $\psi$
and $\theta:=\psi^{-1}\circ \phi$ be isomorphisms
of topological vector spaces.
Since $\phi$ and $\psi$ are isomorphisms
of vector spaces,
we can always identify $(E_1\tensor E_2)\tensor E_3$,
$E_1\tensor (E_2\tensor E_3)$
and $E_1\tensor E_2\tensor E_3$
as vector spaces for simplicity of notation;
only the topologies may differ.
\begin{prop}\label{criter}
In the preceding situation, we have:
\begin{itemize}
\item[\n (a)]
If $E_3$ is locally bounded, then $\phi$
is an isomorphism of topological vector spaces.
\item[\n (b)]
If $E_1$ and $E_3$ are locally bounded,
then the natural isomorphism of vector spaces
\[
\theta\!: (E_1\tensor E_2)\tensor E_3\to
E_1\tensor (E_2\tensor E_3)
\]
taking
$(x\tensor y)\tensor z$ to $x\tensor (y\tensor z)$
is an isomorphism of topological vector spaces.
\end{itemize}
\end{prop}
\begin{proof}
(a) Let $B$ be a bounded, balanced zero-neighbourhood
in~$E_3$. Then $\{tB\!: t\in \R^\times\}$ is a basis of
zero-neighbourhoods in~$E_3$.
Hence
a basis of zero-neighbourhoods of $E_1\tensor E_2\tensor E_3$
is given by the sets of the form
\[
\sum_{n\in \N} U_n\tensor V_n\tensor (t_nB)
=
\sum_{n\in \N} (t_nU_n)\tensor V_n\tensor B\,,
\]
where $U_n$ and $V_n$ are open zero-neighbourhoods in $E_1$
and $E_2$, respectively, and $t_n\in \R^\times$.
Replacing $U_n$ with $t_nU_n$,
we see that it suffices to take $t_n=1$
for all $n\in \N$ here:
the sets
\[
W:=\sum_{n\in \N} U_n\tensor V_n\tensor B
\]
form a basis of zero-neighbourhoods. Suppose such a $W$
is given. We choose a bijection
$\lambda\!: \N\to \N^2$
and define
$P_{\lambda(n)}:=U_n$ and $Q_{\lambda(n)}:=V_n$
for $n\in \N$.
Then
\[
W
%=\sum_{n\in \N} U_n\tensor V_n\tensor B
=\sum_{n\in \N} P_{\lambda(n)}\tensor Q_{\lambda(n)}\tensor B
=\sum_{(j,k)\in \N^2} P_{(j,k)}\tensor Q_{(j,k)}\tensor B
\supseteq
\sum_{j \in \N} \Bigl(\sum_{k\in \N} P_{(j,k)}\tensor Q_{(j,k)}\Bigr)
\tensor B\,,
\]
where the right hand side is a zero-neighbourhood
in $(E_1\tensor E_2)\tensor E_3$. Thus $\phi$ is open
and thus $\phi$ is an isomorphism of topological vector spaces.\vspace{2mm}

(b) If also $E_1$ is locally bounded,
we see in the same way that
$\psi$ is an isomorphism of topological vector spaces,
whence so is
$\theta=\psi\circ \phi^{-1}$.
\end{proof}
\section{Examples where associativity fails}
We show that
tensoring in the category of topological vector spaces
is not associative.
\begin{thm}\label{main}
Let $E:=\R^\N$ with the product topology.
Then the following holds:
\begin{itemize}
\item[\n (a)]
The canonical isomorphism of vector spaces
\[
\theta\!: (E\tensor E)\tensor E\to E\tensor (E\tensor E)
\]
is not continuous $($and hence not an isomorphism
of topological vector spaces$)$.
\item[\n (b)]
The canonical isomorphisms of vector spaces
\[
\phi\!:E\tensor E\tensor E\to (E\tensor E)\tensor E\qquad
\mbox{and}\qquad
\psi\!: E\tensor E\tensor E\to E\tensor (E\tensor E)
\]
are continuous but not open.
\end{itemize}
\end{thm}
\begin{proof}
(a)
Given $\ve>0$ and $n\in \N$, the set
$U_{n,\ve}:= \{z\in \R\!: |z|<\ve\}^n\times \R^{\{n+1,n+2,\ldots\}}$
is a zero-neighbourhood in~$E$, whence
\[
U:=\sum_{n\in \N}\, U_{2n, 2^{-n}}\tensor
\Bigl(\sum_{k\in \N} \,\bigl( U_{2n,2^{-k}}
\tensor U_{2n, 1}\bigl) \Bigr)
\]
is a zero-neighbourhood in $E\tensor (E\tensor E)$.
Then $U$ is not a zero-neighbourhood
in $(E\tensor E)\tensor E$.
In fact, otherwise we find
zero-neighbourhoods $A_{n,k}$, $B_{n,k}$
and $C_n$ in~$E$ such that
\[
V:=\sum_{n\in \N}\Bigl(\sum_{k\in \N} A_{n,k}\tensor B_{n,k}\Bigr)
\tensor C_n\sub U\,.
\]
There is $N\in 2 \N$ such that $\R\, e_N\sub
C_1$, where $e_N:=\delta_{N,\sbull}\!: \N\to \{0,1\}\sub \R$
is defined using Kronecker's delta.
The sets $A_{1,k}$ and $B_{1,k}$ being absorbing,
we then have
\begin{equation}\label{allin}
\Bigl(\sum_{k\in \N} E\tensor E\Bigr) \tensor \R \, e_N \sub V\sub U\,.
\end{equation}
Let $p\!: E\tensor E\tensor E\to \R^N\tensor \R^N$
be the linear map uniquely determined by
\[
p(x\tensor y\tensor z)=z(N)\cdot
\bigl(x|_{\{1,\ldots, N\}}\bigr)
\tensor \bigl(y|_{\{1,\ldots, N\}}\bigr)\,.
\]
We now identify $\R^N$ with
$\R^N\times \{0\}\sub \R^\N$.
Given $v\in \R^N\tensor \R^N$
and $R\in \R^\times$,
we have $(Rv) \tensor e_N\in U$ by (\ref{allin}),
whence $(Rv)\tensor e_N=\sum_{k,n\in \N}
x_n\tensor y_{n,k}\tensor z_{n,k}$
for suitable
$x_n\in U_{2n,2^{-n}}$,
$y_{n,k}\in U_{2n,2^{-k}}$,
and $z_{n,k}\in U_{2n,1}$
(depending on~$R$),
almost all of which are zero.
Hence
\[
v=p(v\tensor e_N)=\sum_{k,n\in \N} {\textstyle \frac{c_{n,k}}{R}}
\cdot u_n\tensor v_{n,k}
=\sum_{n\in \N} u_n\tensor
\underbrace{\Bigl(
\sum_{k\in \N}{\textstyle \frac{c_{n,k}}{R}}
\cdot v_{n,k}\Bigr)}_{=:w_n}
\]
with $u_n:=x_n|_{\{1,\ldots, N\}}$,
$v_{n,k}:=y_{n,k}|_{\{1,\ldots, N\}}$ and $c_{n,k}:=z_{n,k}(N)$.
Thus
\begin{equation}\label{zuwenige}
v= \sum_{n=1}^{N/2} u_n\tensor w_n +\rho
\end{equation}
with $\rho:=\sum_{n>N/2} u_n\tensor w_n
=\sum_{i,j=1}^N \Bigl( \sum_{n>N/2}\sum_{k\in \N}
{\textstyle \frac{c_{n,k}}{R}} u_n(i)v_{n,k}(j)\Bigr)
e_i\tensor e_j$. Note that,
as $2n>N\geq i,j$,
we have $|u_n(i)|<2^{-n}$, $|v_{n,k}|<2^{-k}$
and $|c_{n,k}|<1$. Thus
\[
\sum_{n>N/2}\sum_{k\in \N}
\bigl|{\textstyle \frac{c_{n,k}}{R}} u_n(i)v_{n,k}(j)\bigr|
\leq
{\textstyle \frac{1}{|R|}}
\sum_{n,k=1}^\infty 2^{-n}2^{-k}
={\textstyle \frac{1}{|R|}}\,,
\]
which can be made arbitrarily small
for large $|R|$. Thus (\ref{zuwenige})
shows that $v$ is contained in the closure
of the set~$S$ of $N/2$-fold sums of elementary
tensors, with respect to the canonical Hausdorff
vector topology ($\isom \R^{N^2}$) on $\R^N\tensor \R^N$.
Hence $S$ is dense in $\R^N\tensor \R^N$.
However, under the usual
isomorphism $\R^N\tensor \R^N\isom M_N(\R)$
the set~$S$ corresponds to the set of matrices
of rank~$\leq N/2$, which is not dense in $M_N(\R)$
as its closure does not meet the open set $\GL_N(\R)$
of invertible matrices. We have reached a contradiction.\vspace{2mm}

(b) If one of $\phi$ and $\psi$
was an isomorphism of topological vector spaces,
then also the other,
because $\psi=\alpha\circ \phi\circ \beta$,
where the linear maps
$\beta\!: E\tensor E\tensor E\to E\tensor E\tensor E$
and $\alpha\!: (E\tensor E)\tensor E\to E\tensor (E\tensor E)$
determined by $\beta(x\tensor y\tensor z):=y\tensor z\tensor x$
and $\alpha((x\tensor y)\tensor z)=z\tensor (x\tensor y)$
are isomorphisms of topological vector spaces.
Thus $\theta=\psi\circ \phi^{-1}$ would be an
isomorphism of topological vector spaces,
contradicting\,(a).
\end{proof}
\begin{rem}\label{conseq}
The spaces
$E\tensor E$, $(E\tensor E)\tensor E$
and $E\tensor (E\tensor E)$
are also the respective tensor products
in the category of {\em Hausdorff\/} real topological vector
spaces (Proposition~\ref{propo}\,(d)),
whence
Theorem~\ref{main}\,(a) entails that tensoring is
not associative in this (more interesting) category.
Nor is tensoring associative in the category
of complete Hausdorff real topological vector spaces,
since $E\tensor E$, $(E\tensor E)\tensor E$
and $E\tensor (E\tensor E)$
are complete (see \cite{Tu2}).
Likewise, considering
$E\tensor E$, $(E\tensor E)\tensor E$
and $E\tensor (E\tensor E)$
as tensor products in the category of
Hausdorff abelian topological groups
(Proposition~\ref{propo}\,(e)),
we deduce that tensoring is not associative
in the category of Hausdorff abelian topological
groups.
\end{rem}
\begin{rem}
Note that the definition of tensor products
and all of the results obtained (except for
Proposition~\ref{propo}\,(d),\,(e)
and Remark~\ref{conseq})
remain meaningful and correct when $\R$ is replaced with
an arbitrary complete valued field~$\K$ (for instance,
the field of $p$-adic numbers).
This observation is of interest
in connection with
\cite{BGN},
where a framework of
differential calculus over non-discrete
topological fields is described.
In this context,
it is desirable to know precisely which constructions
of topological algebra work over general topological fields
(or at least complete valued fields),
in contrast to those which
depend on specific properties of the real number field,
or on local convexity.
\end{rem}
{\em Open problems}. Does Proposition~\ref{propo}\,(d) carry over
to Hausdorff topological vector spaces over complete valued
fields~$\K$ other than $\R$ (or $\C$)\,? Are tensor products
of complete Hausdorff topological $\K$-vector spaces
in the class of Hausdorff topological $\K$-vector spaces
always complete, as in in the real case (discussed in \cite{Tu2})\,?\\[3mm]
{\em Acknowledgements.} The author wishes to thank K.-H Neeb
for a suggestion which lead to
Proposition~\ref{propo}\,(e).
{\footnotesize
{\bf Helge Gl\"{o}ckner}, TU~Darmstadt, FB~Mathematik~AG~5,
Schlossgartenstr.\,7, 64289 Darmstadt, Germany.\\
E-Mail: gloeckner@mathematik.tu-darmstadt.de}
\end{document}